\documentclass[12pt, amssymb]{article}

\usepackage{latexsym}   
\usepackage{amssymb}    
\usepackage{amsmath}    
\usepackage{amsthm}     

\newtheorem{theorem}{Theorem}[section]
\newtheorem{lemma}[theorem]{Lemma}

\def\X{{\cal X}}
\def\Y{{\cal Y}}
\def\O{{\cal O}}
\def\D{{\cal D}}

\def\Z{{\mathbb Z}}

\def\Q{{\mathbb Q}}

\def\P{{\mathbb P}}

\begin{document}
\title{Invariants Associated to Orthogonal $\epsilon$-constants}
\author{Darren Glass}
\date{}
\maketitle

\begin{abstract} In this paper we use the theory of $\epsilon$-constants
associated to tame finite group actions on arithmetic surfaces to define a
Brauer group invariant $\mu(\X,G,V)$ associated to certain symplectic motives
of weight one.  We then discuss the relationship between this invariant and
$w_2(\pi)$, the Galois theoretic invariant associated to tame covers of
surfaces. \end{abstract}

\section{Introduction and Background}
\label{section:intro}
\setcounter{equation}{0}

Brauer group invariants associated to motives have been long studied by
mathematicians.  In the case of orthogonal motives of even weight, these
invariants were first examined by Frohlich in \cite{F}.  Deligne used these
Brauer group invariants and their relationship with certain
$\epsilon$-constants in order to give a proof of the Frohlich-Queyrut Theorem
in \cite{D2} by showing that certain global orthogonal root numbers were one by
interpreting the associated local orthogonal root numbers as Stiefel Whitney
classes and then using the local root numbers to define an element of order two
in the Brauer group of $\Q$.  Similar results and methods have been done by
Saito (in \cite{tS2}, for example) and others.  In this paper, we will define a
Brauer group invariant associated to certain motives which are symplectic and
have weight one.  

In order to construct the relevant motives we first define $\X$ to be an
arithmetic surface of dimension $2$ which is flat, regular, and projective over
$\Z$.  Throughout this paper we will assume that $f:\X \rightarrow Spec(\Z)$ is
the structure morphism.  Define $G$ to be a finite group which acts tamely on
$\X$.  In other words, for each closed point $x \in X$, the order of the
inertia group of $x$ is relatively prime to the residue characteristic of $x$. 
Let $\Y$ be the quotient scheme $\X / G$, which we will assume is regular, and
that for all finite places $v$  the fiber $\Y_v = (\X_v)/G = \Y
\otimes_\Z(\Z/p(v))$ has normal crossings and smooth irreducible components
with multiplicities relatively prime to the residue characteristic of $v$.
Finally, let $V$ be a representation of $G$ over $\overline{\Q}$

We wish to define a class $\mu(\X,G,V)$ in the global Brauer group
$H^2(\Q,\Z/2\Z)$ which will depend only on $\X, G,$ and $V$.  In particular,
$\mu(\X,G,V)$ will be the element whose local invariant at the place $v$ is
given by the sign of $\epsilon_v(D',V)$, where $D'$ is an appropriately chosen
horizontal canonical divisor on $\Y$.  In other words, $O_\Y(D' + \Y_T)$ is
isomorphic to $\omega_{\Y/\Z}(\Y_S^{red})$ where $\Y_T$ and $\Y_S$ are unions
of vertical fibers of $\Y$.

In addition to being a canonical divisor, we must impose an additional
condition on our choice of $D'$ in order for $\mu(\X,G,V)$ to be well-defined. 
To do this we recall that Proposition $3.10$ of \cite{G} showed that if $F'$
and $G'$ are components of $\Y_v^{red}$ and $D'$ is a canonical divisor in the
above sense, then there is a canonical isomorphism between $\O_{F'}(D' \cap
F')$ and $\omega_{F'}(F' \cap G')$.  This isomorphism maps the global section
$1 \in \Gamma(\O_{F'}(D' \cap F'))$ to an element $\gamma \in
\Gamma(\omega_{F'}(F' \cap G'))$ such that $\gamma$ has a simple pole at all $x
\in F' \cap G'$.  Define $a_x$ to be the residue of $\gamma$ at these points
$x$.  We wish to choose $D'$ so that all of the residue terms $a_x$ are equal
to $1$ and we will call such a choice of $D'$ a nice divisor.   Nice divisors
exist due to a moving lemma which we will make explicit in Section 2.  We can
now state the main theorem.

\begin{theorem} 

\label{prop:brauer}   In the case where $V$ is an orthogonal
representation of dimension equal to zero and trivial determinant and
where $D'$ is a nice divisor, the local constant $\epsilon_v(D',V)$ is
independent of the choice of $D'$.  In particular, the element
$\mu(\X,G,V)$ in the global Brauer group $H^2(\Q,\Z/2\Z)$ whose local
invariant at the place $v$ is given by the sign of $\epsilon_v(D',V)$
is well-defined.

\end{theorem}

Section two of this paper will prove theorem \ref{prop:brauer}.  Section three
recalls the definition of the Galois theoretic invariant $w_2(\pi)$ defined by
Cassou-Nogues, Erez and Taylor in \cite{CNET} which comes from a situation
similar to the one we are working in and proves a connection between the two
invariants.  Finally, in section four we prove the moving lemma which we used
in section two in order to show that nice divisors always exist.

The author would like to thank Ted Chinburg for his insightful comments and
suggestions, in particular on the proof in section four.

\section{Definition of $\mu(\X,G,V)$}

In \cite{G}, the author proves the following result about orthogonal
$\epsilon$-constants associated to tame finite group actions on surfaces (We
refer the reader to \cite{D1}, \cite{CEPT1}, and \cite{G} for all relevant
definitions.): 

\begin{theorem}
\label{thm:thesis}

If $V$ is an orthogonal virtual representation of degree zero and trivial
determinant then for all finite places $v$ of $\Q$ we have the following
formula: 
\begin{eqnarray*}
\label{eqn:prodform}
\frac{\epsilon(\Y_v,V)}{\epsilon(D'_v,V)}
=\prod_i det(V^{I_i})(\delta_{v,C_i})\prod_{z \in Z} \epsilon_{0,z}(C_{i_1},V^{I_{i_1}})\epsilon_{0,z}(C_{i_2},V^{I_{i_2}})
\epsilon(z,V)
\end{eqnarray*}
where $Z$ is the set of crossing points of the components of the fiber $\Y_v$
and $D'$ is chosen to be a canonical horizontal divisor on $\Y$.
\end{theorem}

By canonical we mean that $D'$ is a horizontal divisor on $\Y$ such that $D'+ \Y_T = K_{\Y} + \Y_S^{red}$, where $\Y_S^{red}$ is the sum of the reductions of the fibers of $\Y$ at the places in the set $S$ of bad primes, $\Y_T$ is the sum of the (necessarily reduced) fibers of $\Y$ over the places in $T$.  We also wish to choose $D'$ so that it intersects the non-smooth fibers $\Y_v$ of $\Y$ transversally at smooth points on the reduction of $\Y_v$. 

Given any such horizontal divisor $D'$ it is possible to find a nice divisor which is close to it due to the following moving lemma, which is proven in section \ref{section:moving}.  In particular, this shows that nice divisors $D'$ always exist and therefore that our class $\mu(\X,G,V)$ will be well-defined.

\begin{lemma} \label{lemma:moving} There exists a meromorphic function
$h$ on $\Y_v^{red}$ such that the divisor of $h$ intersects the special
fibers $\Y_v^{red}$ transversally at smooth points away from $D_v'$ and
such that $h$ takes on prescribed values at the singular points of
$\Y_v^{red}$.  In particular, given a horizontal divisor $D'$ as in the
previous section, the divisor $D' + div(h)$ will have residue maps equal
to one at the crossing points of components of $\Y_v^{red}$.  \end{lemma}

\noindent {\bf Proof of theorem \ref{prop:brauer}:}\enspace   
 For any fixed place $v$ of $\Q$ and any component $C_i$ of $\Y_v^{red}$, it follows from the definition that $\delta_{v,C_i}$ is the class which corresponds to the element $\delta=(\oplus 0) \oplus (\oplus a_x) \in (\oplus_{x \in C_i - Z}\Z) \oplus (\oplus_{x \in C_i \cap Z}K^*/U^1_x)$.  This term is independent of the choice of nice divisor and it follows that the right hand side of the equation in Theorem \ref{thm:thesis} is as well. Furthermore, it is clear that $\epsilon(\Y_v,V)$ is independent of our choice of $D'$.  Thus, it follows from the theorem that $\epsilon_v(D',V)$ is independent of the
choice of $D'$. Next, we note that $\epsilon_{v,0}(D',V)$ must also be independent of our choice of $D'$ (by Lemma 3.3 of \cite{G}, for example). 
Therefore it must be the case that $\epsilon_v(D',V)=
\epsilon_{v,0}(D',V)\epsilon(D_v',V)$ is independent of the choice of
$D'$. The product of all of the $\epsilon_v(D',V)$ is equal to
$\epsilon(D',V)$, which must be equal to one from the
theorem of Fr\"ohlich and Queyrut \cite{FQ}.  This tells us that we can define an
element $\mu(\X,G,V)$ in $H^2(\Q,\Z/2\Z)$ by setting the local component
at the prime $v$ to be equal to the sign of $\epsilon_v(D',V)$.  $\blacksquare$

\section{The connection to $w_2(\pi)$}

Let $\pi: \X \rightarrow \Y$ be a tamely ramified cover of degree $n$,
where $\X$ and $\Y$ are regular schemes and $\Y$ is connected.
Furthermore, we must make the technical assumption that the
ramification indices are all odd. Cassou-Nogues, Erez, and Taylor use
Grothendeick's equivariant cohomology theory to define an invariant
$w_2(\X/\Y) = w_2(\pi) \in H^2(\Y_{et},\Z/2\Z)$ associated to this
situation. Their definition generalizes to define classes $w_i(\pi)$
which lie in $H^i(\Y_{et},\Z/2\Z)$ for all positive integers $i$, but
in this thesis we will only be interested in $w_2$. These terms are
generalized Stiefel-Whitney classes, and are obtained by pulling back
the universal Hasse-Witt classes defined by Jardine using classifying
maps related to a quadratic form $E$. The precise definition of $E$
uses the existence of a locally free sheaf $\D_{\X/\Y}^{-1/2}$ whose
square is the inverse different of the covering $\X/\Y$. In this section, we will consider the relationship between the class $\mu(\X,G,V)$ lying in $H^2(\Q,\Z/2\Z)$ which we defined in theorem \ref{prop:brauer} and $w_2(\pi)$.

In \cite{CNET} the following equality in
$H^2(\Y_{et},\Z/2\Z)$, which is an analog of a theorem of Serre, is proved:
$$w_2(\pi_*(\D_{\X/\Y}^{-1/2},Tr_{\X/\Y})) = w_2(\pi) + (2) \cup
(d_{\X/\Y}) + \rho(\X/\Y)$$

\noindent where $\rho(\X/\Y)$ is defined by the ramification of
$\X/\Y$, $d_{\X/\Y}$ is the function field discriminant, and the left
hand side of the equation is the Hasse-Witt invariant associated to the
square-root of the inverse different bundle.  Note that if we look at
the one-dimensional version of this formula the middle term on the
right hand side becomes trivial.  Therefore, in the case of \'etale
covers of curves the formula reduces to $$w_2(\pi) = w_2(\pi_*(\D_{\X/\Y}^{-1/2},Tr_{\X/\Y})) = w_2(E)$$

\noindent where $w_2(E)$ is the second Hasse-Witt invariant associated
to the square root of the inverse different, as described in detail in
\cite{CNET}.

Let $D'$ be a choice of a canonical divisor on $\Y$ in the sense of the
previous sections, and let $i:D' \hookrightarrow \Y$ be the natural
inclusion.   An \'etale covering of $\Y$ naturally restricts to give an
\'etale covering of $D'$.  We now have the following natural maps

$$i^*:H^2(\Y_{et},\Z/2\Z)\rightarrow H^2(D'_{et},\Z/2\Z)$$

$$res: H^2(D'_{et},\Z/2\Z) \rightarrow H^2_{et}(\Q(D'),\Z/2\Z) =
H^2_{gal}(\overline{\Q}/\Q(D'),\Z/2\Z)$$

$$cor:H^2_{gal}(\overline{\Q}/\Q(D'),\Z/2\Z)\rightarrow
H^2(\Q,\Z/2\Z)$$

\noindent where the latter two maps are restriction and corestriction
in the sense of Serre (for details see Chapter VII of \cite{Se}). 
Composing these maps gives a natural map

$$H^2(\Y_{et},\Z/2\Z)\rightarrow H^2(\Q,\Z/2\Z)$$  

We denote the image of the class $w_2(\pi) \in H^2(\Y_{et},\Z/2\Z)$
under this map by $\tilde{w}_2(\pi) \in H^2(\Q,\Z/2\Z)$. At first
glance it appears as though this element may depend on our choice of
canonical divisor $D'$.  However, we want to show that it does not
depend on this choice and furthermore that the element $\tilde{w}_2(\pi)$ is connected in a natural way to the element $\mu(\X,G,V)$. Recall that $\mu(\X,G,V)$ is defined by letting the local invariant at the place $v$ be given by the sign of $\epsilon(D'_v,V)$.  This also appears to depend on the choice of $D'$ but turns out to be independent of choice. 

The first natural observation is that the class $\mu(\X,G,V)$ depends on the choice of a
representation $V$ of $G$ while $\tilde{w}_2(\pi)$ does not.  The natural representation to consider is the regular representation of the group $G$, which we denote by $R$.  In particular, the nicest possible theorem comparing the invariants would say that $\mu(\X,G,R) = \tilde{w}_2(\pi)$.
However,we have only shown that $\mu(\X,G,V)$ is a well-defined class
in the case where $V$ is of dimension zero and of trivial determinant,
neither of which holds for $R$. So instead of setting $V=R$, we
consider the representation $V = R - det(R) - T^{n-1}$, where $det(R)$,
the determinant of the regular representation, is a character whose
order is either one or two, $T$ is the trivial representation and $n$
is the degree of the cover $\X/\Y$.  This choice of $V$ is an
orthogonal representation, and it has trivial determinant and dimension
$0$. We can now prove the following theorem:

\begin{theorem}
\label{thm:cnet-dg}

Assume that we are in the above situation, and in particular that $V =
R - det(R) - T^{n-1}$.  Let $\Y_1/\Y$ be either the trivial cover or
the subcover of $\X/\Y$ of degree 2, depending on whether $det(R)$ is of order $1$ or $2$ respectively.  Then as classes in $H^2(\Q,\Z/2\Z)$, we have the equality  $$\mu(\X,G,V) =
\tilde{w}_2(\X/\Y) - \tilde{w}_2(\Y_1/\Y) - (n - 1)\tilde{w}_2(\Y/\Y)$$

\end{theorem}

The proof of this theorem relies on the interpretation of each side of
the equation as a Stiefel-Whitney class.  In particular, Cassou-Nogues,
Erez, and Taylor show that the element $w_2(\pi)$ is the Hasse-Witt
invariant associated to the full covering of surfaces.  Thus, when we
restrict the class to the one-dimensional divisor $D'$ we see that the
element $i^*(w_2(\pi)) \in H^2(D'_{et},\Z/2\Z)$ is equal to the
Stiefel-Whitney class associated to the form $E' =
(\D_{D/D'}^{-1/2},Tr_{D/D'})$ on the canonical divisor $D'$ of $\Y$. 
This follows as a generalization to \'etale cohomology of results of
Fr\"ohlich in \cite{F}, which allow us to associate the class
$i^*(w_2(\pi))$ to $G$-extensions of the ring of integers of the
residue field of the generic point of $D'$.

Next we make use of the results of Deligne which allow us to interpret
local Stiefel-Whitney classes in terms of local root numbers.  In
particular, the following lemma is shown in {\cite{D2}}:

\begin{lemma}
\label{lemma:deligne}

Let $d=1$, so that the fibers $\X_v$ and $\Y_v$ are all one dimensional
schemes.  Furthermore, let $V$ be an orthogonal virtual representation
of dimension zero and trivial determinant.  Under these hypotheses, the
local root number $W(V_v) = sign(\epsilon_v(\Y,V))$ is equal to $exp(2
\pi i cl(sw_v))$, where $sw_v$ is the local Stiefel-Whitney class, and
$cl(sw_v) \in \{0, 1/2\} \subset \Q/\Z$.

\end{lemma}

In other words, in characteristic not equal to two, the sign of the
$\epsilon$-constants $\epsilon_v(D',V)$ of the representation on the
one dimensional horizontal divisor $D'$ are determined by whether or
not the classes $w_2(\pi)$ are trivial in the Brauer group, and
$\epsilon_v(D',V)$ is automatically positive when $v=2$. It turns out
that these are exactly the terms that are coming up in the computation
of the class of Cassou-Nogues, Erez, and Taylor.  

In particular, $\epsilon_v(D',V)=\epsilon_v(D',R)\epsilon_v(D',det(R))$
is the same as the local Hasse-Witt invariants.  However, we are
working with \'etale covers of curves and so from the results of
\cite{CNET} discussed above, these Hasse-Witt invariants are simply the
images of the appropriate classes $w_2(\pi)$.  This proves Theorem
\ref{thm:cnet-dg}. $\blacksquare$









\section{Proof of Lemma \ref{lemma:moving}}
\label{section:moving}

The proof of Lemma \ref{lemma:moving} involves a generalized version of
Bertini's Theorem.  For now, let us assume that $X$ is a smooth curve
defined over an infinite field $k$ and let us choose a finite set of
points $p_1,\ldots,p_m$ on $X$.  We define the divisor $p = \sum_i p_i$.
Furthermore, let us choose constants $c_i$ which lie in the residue
field $k(p_i)$ of the points $p_i$.   Finally, let us choose  $\Lambda$
to be an effective very ample divisor on $X$ of large degree which is
supported off of $p$.  We look at the group of global sections $H^0(X,
\O_X(\Lambda))$ and let $f_0,\ldots,f_t$ be a basis of this group.  This
basis gives us a projective embedding from $X$ into $\P_k^t$ whose projective
coordinates we will write as $x_0,\ldots,x_t$.  

We wish to prove that there exist linear forms $l_0$ and $l_1$ in the
variables $x_i$ such that the following properties hold:

\begin{enumerate}

\item[(1)] If $H_i$ is the hyperplane defined by
$l_i=0$ in $\P_k^t$ then $H_i \cap X$ is a finite set of closed points
which is regular and disjoint from $\{p_1,\ldots,p_m\}$.  Furthermore,
we wish to choose the $l_i$ so that $H_1 \cap H_2 \cap X$ is empty.

\item[(2)] It follows from (1) that the function $l_1/l_0 |_X$ is in
$\O_{X,p_i}$ for each $i$.  However, we further wish to specify that the
image of $l_1/l_0$ in $k(p_i)$  is the given constant $c_i$. 
\end{enumerate}

The classical version of Bertini's Theorem (Theorem $II.8.18$ of \cite{H})
tells us that there exist linear forms $l_0$ so that $H_0$ satisfies
condition (1).  We now fix one choice of such an $l_0$, and we will
attempt to construct an $l_1$ so that the pair satisfies properties (1)
and (2). We begin by looking at the set $V$ consisting of all linear
forms such that $\{l_0,l_1\}$ satisfy condition (2).  In other words,

$$V=\{l = a_0x_0+\ldots+a_tx_t \mid \forall j, \frac{l}{l_0}|_X(p_j)=c_j \in k(p_j)\}$$ 

This $V$ will be an affine space over $k$. Furthermore, because we
chose the divisor $\Lambda$ to have high degree it follows from a
Riemann-Roch argument that $V$ is of codimension $m$ inside of $H^0(X,
\O_X(\Lambda))$. 

For each point $x \in X$, we now define a set $V_x \subseteq V$ which
consists of all linear forms $l \in V$ so that the hyperplane defined by
${l=0}$ has contact order $>1$ at $x$. In other words, $V_x$ will
consist of those linear forms who do not intersect $X$ nicely at the
point $x$.  We can again use the Riemann-Roch
theorem to show that for almost all choices of $x$, we get that the
dimension of $V_x$ is equal to $dim$ $V - 2$. 

Let $U = X - \{p_1,\ldots,p_m\}$ so that $U$ is an affine curve, and 
define $T \subseteq U \times V$ to be the set of all pairs $(x,l)$ such
that $x \in U$ and $l \in V_x$.  We have seen that the projection map
$\pi:T \rightarrow U$ is surjective and for almost all $x \in U$ (in
particular for those points such that $k(x)=k$), we see that the fiber
$\pi^{-1}(x)$ is an affine space whose dimension is equal to $dim$
$V-2$.  In particular, this shows that $T$ is irreducible and that the
dimension of $T$ is equal to $dim$ $V-1$.  But this shows us that the
natural projection map $\gamma:T \rightarrow V$ must not be surjective.  

In particular, we can choose some element $l_1 \in V$ which is not
in the image of $\gamma$.  In particular, the hyperplane $H$ defined by
$\{l_1 =0\}$ is such that $H \cap U$ is regular and, since $l_1 \in V$,
we know that $l_1/l_0(p_i) = c_i \in k(p_i)$, and thus that $l_1$ and
$l_o$ satisfy conditions (1) and (2) above.

So far, we have only made the argument for the case where $X$ is a
smooth curve.  However, as long as $X$ is a reduced curve with smooth
irreducible components which have normal crossings, then the same
argument will hold as long as we include these crossing points in the
set of $\{p_i\}$.  Instead of using the normal Riemann-Roch theorem we
will now use the version for singular curves described on p.298 of
\cite{H}.

In order to prove lemma \ref{lemma:moving} we apply this generalized
version of Bertini's theorem to $X=\Y_v^{red}$.  Specifically,
choose the set of points $\{p_1,\ldots,p_m\}$ to include the crossing
points of components of $\Y_v^{red}$ as well as the points in $D' \cap
\Y_v^{red}$.  The above argument then allows us to find a meromorphic
function $h$ where we can specify the values of the function $h =
l_1/l_0$ at the crossing points of $\Y_v^{red}$ so that the residues
that come up  when we consider $D'' = D' + div(h)$ are all equal to
one and $D''$ intersects $\Y_v^{red}$ in the way we want.  $\blacksquare$

\end{document}